\documentclass{jgcc}
\pdfoutput=1
\usepackage{lastpage,orcidlink,hyperref}

\jgccdoi{18}{2}{2}{16903}

\jgccheading{}{\pageref{LastPage}}{}{}{Nov.~11,~2025}{Jun.~17,~2026}{}


\DeclareMathOperator{\ar}{ar}
\DeclareMathOperator{\Tm}{Tm}
\DeclareMathOperator{\supp}{supp}

\newcommand{\B}{\{0,1\}}
\newcommand{\cl}{\mathfrak}
\newcommand{\bb}{\mathsf}
\newcommand{\cs}{,\,}
\newcommand{\st}{\,|\,}
\newcommand{\rep}[2]{[#1]_{#2}}
\newcommand{\Add}[1]{#1|_\Gamma}
\newcommand{\alg}[1]{\langle#1\rangle}
\newcommand{\gp}[1]{\langle#1\rangle_\Gamma}
\newcommand{\GP}[1]{\left\langle#1\right\rangle_\Gamma}
\newcommand{\acl}[1]{\lvert#1\rvert}
\newcommand{\N}{\mathbb N}
\newcommand{\Z}{\mathbb Z}
\newcommand{\R}{\mathbb R}
\newcommand{\rv}{\mathbf}
\newcommand{\re}{\mathrm}
\newcommand{\Prob}[1]{\Pr[#1]}
\newcommand{\Forall}[1]{\forall\,#1\ }
\newcommand{\Exists}[1]{\exists\,#1\ }

\newcommand{\eg}{\emph{e.g.}}
\newcommand{\cf}{\emph{cf.}}

\begin{document}
\title[Polynomial-Time Algorithms]{Polynomial-Time
Algorithms for Black-Box Distributive Expanded Groups}

\author[M.~Anokhin]{Mikhail
Anokhin\,\orcidlink{0000-0002-3960-3867}}

\address{Information Security Center, Faculty of
Computational Mathematics and Cybernetics, Lomonosov Moscow
State University, Moscow, Russia}

\email{anokhin@mccme.ru}


\begin{abstract}
Let $\Omega$ be a finite set of finitary operation symbols.
An $\Omega$-expanded group is a group (written additively
and called the additive group of the $\Omega$-expanded
group) with an $\Omega$-algebra structure. We use the
black-box model of computation in $\Omega$-expanded groups.
In this model, elements of a finite $\Omega$-expanded group
$H$ are represented (not necessarily uniquely) by bit
strings of the same length, say,~$n$. Given representations
of elements of $H$, equality testing and the fundamental
operations of $H$ are performed by an oracle.

Assume that $H$ is distributive, i.e., all its fundamental
operations associated with nonnullary operation symbols in
$\Omega$ are distributive over addition. Suppose
$s=(s_1,\dots,s_m)$ is a generating system of~$H$. In this
paper, we present probabilistic polynomial-time black-box
$\Omega$-expanded group algorithms for the following
problems: (i)~given $(1^n,s)$, construct a generating system
of the additive group of $H$, (ii)~given
$(1^n,s,(t_1,\dots,t_k))$ with $t_1,\dots,t_k\in H$, find a
generating system of the additive group of the ideal in $H$
generated by $\{t_1,\dots,t_k\}$, and (iii)~given $(1^n,s)$,
decide whether $H\in\cl V$, where $\cl V$ is an arbitrary
finitely based variety of distributive $\Omega$-expanded
groups with nilpotent additive groups. The error probability
of these algorithms is exponentially small in~$n$. In
particular, this can be applied to groups, rings,
$R$-modules, and $R$-algebras, where $R$ is a fixed finitely
generated commutative associative ring with~$1$. Rings and
$R$-algebras may be here with or without $1$, where $1$ is
considered as a nullary fundamental operation.
\end{abstract}

\maketitle

\section{Introduction}

In~\cite{BS84}, Babai and Szemer\'edi introduced a model of
computation in finite groups, called the black-box group
model. In that model, elements of a finite group $G$ are
represented (not necessarily uniquely) by bit strings of the
same length. Given representations of elements of $G$,
equality testing and the group operations in $G$ are
performed by an oracle. When an algorithm performs a
computation in $G$, it works with the above-mentioned
representations of elements of $G$ and has access to an
arbitrary oracle performing the above tasks. Also, this
algorithm takes as input a tuple $(r_1,\dots,r_m)$, where
$r_i$ is a representation of $s_i\in G$ for every
$i\in\{1,\dots,m\}$; it is assumed that the system
$(s_1,\dots,s_m)$ generates the group~$G$.

We note that in~\cite{BS84}, testing equality of elements of
$G$ (more precisely, the equivalent task of testing equality
to $1$ in $G$) is understood in a special sense. However, we
understand testing equality of elements of $G$ in the
standard sense, i.e., given representations of arbitrary
elements $g_1,g_2\in G$, decide whether $g_1=g_2$.

Let $\Sigma$ be a fixed finite set of finitary operation
symbols. Then the black-box group model can be naturally
generalized to $\Sigma$-algebras (see
Subsection~\ref{ss:bbunivalgmodel}). Namely, we define a
black-box $\Sigma$-algebra as a pair $(H,\rho)$, where $H$
is a finite $\Sigma$-algebra and $\rho$ is a function from a
subset of $\B^n$ onto $H$ for some $n\in\N$. (Throughout
this paper, $\N=\{0,1,\dotsc\}$.) Suppose $\bb H=(H,\rho)$
is a black-box $\Sigma$-algebra. Any algebraic statement
concerning $\bb H$ should be understood as concerning~$H$.
If $H\ne\emptyset$, then there exists a unique $n\in\N$ such
that $\rho$ is defined on a subset of $\B^n$; this
nonnegative integer $n$ is called the encoding length
of~$\bb H$ (\cf~\cite{BB99}). Otherwise any $n\in\N$ can be
taken as an encoding length of~$\bb H$. For every $h\in H$,
an arbitrary preimage of $h$ under $\rho$ is considered as a
representation of $h$ for computational purposes. For
brevity, we often identify elements of $H$ with their
representations. Loosely speaking, a black-box
$\Sigma$-algebra algorithm is an algorithm that works with
such representations and has access to an arbitrary oracle
performing equality testing and the fundamental operations
during a computation in an arbitrary black-box
$\Sigma$-algebra.

\subsection{Previous work}
\label{ss:prevwork}

There are many works concerning various versions of the
black-box model of computation in algebraic structures. Most
of those works deal with black-box groups. For surveys on
the theory of computation in black-box groups, we refer the
reader to~\cite{BB99} and~\cite[Subsections~2.2--2.4]{CF93};
see also~\cite[Chapter~2]{Ser03}. Black-box prime fields
were studied in~\cite{BL96}. Some polynomial-time algorithms
for associative black-box rings were constructed
in~\cite{ADM06} and~\cite{BP24}. (In this paper, unless
otherwise stated, rings and $R$-algebras are not necessarily
associative, commutative, or with~$1$. Here, $R$ is a
commutative associative ring with~$1$. We do not identify
rings and $\Z$-algebras because they are $\Sigma$-algebras
for different sets~$\Sigma$.)

Note that black-box $\Sigma$-algebra quantum algorithms are
also considered in the literature. See, \eg,~\cite{Wat01,
IMS03} for groups, \cite{ADM06, BP24} for associative rings,
and~\cite{Ano24}. However, in this paper we restrict
ourselves to classical algorithms.

We briefly state some results related to the results of this
paper. Suppose $\bb G=(G,\alpha)$ is a black-box group and
$n$ is the encoding length of~$\bb G$. Babai et
al.~\cite[Theorems~1.5, 2.5, and~3.7]{BCFLS95} proposed
probabilistic polynomial-time black-box group algorithms for
the following problems:
\begin{enumerate}[label=(\roman*)]
\item\label{i:rednumbofelem} Given a generating system
$(s_1,\dots,s_m)$ of $G$, construct a generating system with
$O(n)$ elements for~$G$.

\item\label{i:findgensysfornormcl} Given a generating system
$(s_1,\dots,s_{O(n)})$ for $G$ and a tuple
$(t_1,\dots,t_{O(n)})$ of elements of $G$, find a generating
system with $O(n)$ elements for the normal closure $N$ of
$\{t_1,\dots,t_{O(n)}\}$ in~$G$.
\end{enumerate}
See also~\cite[Subsections~2.3.2, 2.3.3, and~2.4.2]{Ser03}.
Note that if the input is correct, then the algorithm for
problem~\ref{i:rednumbofelem} (resp.,
\ref{i:findgensysfornormcl}) always outputs a system of
elements of $G$ (resp.,~$N$). Therefore, for any constant
$\epsilon$ satisfying $0<\epsilon<1$, we can assume that the
success probability of these algorithms is at
least~$\epsilon$. This can be proved in a standard way by
running the algorithm multiple times and concatenating the
output systems.

Arvind, Das, and Mukhopadhyay~\cite[Theorem~1 and its
proof]{ADM06} constructed a probabilistic polynomial-time
black-box ring algorithm $A$ such that the following holds.
Suppose $\bb R=(R,\gamma)$ is a black-box associative ring
and $s\in R^*$ is a generating system of~$R$. (For any set
$Y$, $Y^*$ denotes the set of all tuples of elements
in~$Y$.) Then $A$ on input $s$ outputs a generating system
of the additive group of $R$ with probability at
least~$1/6$. The algorithm $A$ on input $s$ always outputs a
system of elements of~$R$. Hence for any polynomially
bounded function $\eta\colon\N\to\R_+=\{r\in\R\st r\ge0\}$,
the success probability of $A$ can be amplified to at least
$1-2^{-\eta(n)}$, where $n$ is the encoding length of~$\bb
R$. This can be done in a standard way by running the
algorithm $A$ multiple times and concatenating the output
systems. In this result, $\Sigma$ consists of symbols for
addition, additive inverse, zero, and multiplication.

Let $\cl C$ and $\cl U$ be abstract (i.e., closed under
isomorphism) classes of $\Sigma$-algebras. We say that $\cl
C$ is \emph{G-decidable} in $\cl U$ if there exists a
probabilistic polynomial-time black-box $\Sigma$-algebra
algorithm $A$ such that the following holds. Suppose $\bb
H=(H,\rho)$ is a black-box $\Sigma$-algebra with $H\in\cl
U$, $n$ is an encoding length of $\bb H$ (which is unique if
$H\ne\emptyset$ and arbitrary otherwise), and $s\in H^*$ is
a generating system of~$H$. Then, given $(1^n,s)$ (where
$1^n$ is the unary representation of $n$), $A$ correctly
decides whether $H\in\cl C$ with probability at least~$2/3$.
In Subsection~\ref{ss:bbunivalgmodel}, we explain why we
include $1^n$ in the input to the algorithm $A$ and list
some cases where this is not necessary. The letter G in the
term ``G-decidable'' stands for generating system. We note
that the following conditions are equivalent:
\begin{itemize}
\item The class $\cl C$ is G-decidable in~$\cl U$.

\item There exist a polynomially bounded function
$\xi\colon\N\to\R_+$ and a probabilistic polynomial-time
black-box $\Sigma$-algebra algorithm that decides $\cl C$ in
$\cl U$ as in the above definition of G-decidability, but
with success probability at least $1/2+1/\xi(n)$.

\item For any polynomially bounded function
$\eta\colon\N\to\R_+$, there exists a probabilistic
polynomial-time black-box $\Sigma$-algebra algorithm that
decides $\cl C$ in $\cl U$ as in the above definition of
G-decidability, but with success probability at least
$1-2^{-\eta(n)}$.
\end{itemize}
This equivalence can be proved in a standard way by running
the decision algorithm multiple times, taking the majority
of the outputs, and using the Chernoff bound. In particular,
we obtain an equivalent definition of G-decidability if we
replace the lower bound $2/3$ in the above definition by an
arbitrary constant $\delta$ such that $1/2<\delta<1$.

If we require the algorithm $A$ in the above definition of
G-decidability to be deterministic, then we obtain a
definition of \emph{deterministic G-decidability}; of
course, this deterministic algorithm always gives the
correct answer. These types of decidability (in a slightly
modified form) were defined in~\cite{Ano02} (see
Definition~3.1 in that paper).

Let $k\in\N$. Then $\cl N_k$ denotes the variety of all
nilpotent groups of class at most~$k$. Also, $\cl A$ is the
variety of all abelian groups, i.e., $\cl A=\cl N_1$. If
$\cl V$ and $\cl W$ are varieties of groups, then their
product $\cl V\cl W$ consists of all groups $G$ containing a
normal subgroup $N\in\cl V$ such that $G/N\in\cl W$. It is
well known that multiplication of varieties of groups is
associative. Furthermore, the $k$th power $\cl A^k$ of $\cl
A$ is the variety of all solvable groups of derived length
at most~$k$. See also~\cite[14.64, 14.61, 21.11, 21.51, and
21.52]{Neu67}. In particular, $\cl A^0=\cl N_0$ is the
variety of trivial groups.

Babai et al.~\cite[Corollary~1.6]{BCFLS95} proved that the
classes of all nilpotent and all solvable groups are
G-decidable in the class of all groups. Furthermore, that
corollary implies G-decidability of $\cl N_k$ and $\cl A^k$
in the class of all groups, where $k$ is an arbitrary
nonnegative integer. On the other hand, the classes of all
cyclic and all simple groups are not G-decidable even in the
class of all elementary abelian $p$-groups, where the prime
$p$ depends on the group. This follows from the reasoning of
Babai in~\cite[Section~9]{Bab91}.

In our opinion, the most interesting case of the problem of
deciding $\cl C$ in $\cl U$ is the case where $\cl C$ is
defined within $\cl U$ by a finite set of closed first-order
formulas. This is justified by the fact that in this case,
$\cl C$ satisfies the above definition of deterministic
G-decidability in $\cl U$ with the following modification:
The algorithm $A$ on input $(1^n,(s_1,\dots,s_m))$ runs in
time polynomial in $2^n$ and~$m$. See~\cite[Lemma~3.1
(together with its proof) and Remark~3.4]{Ano02}.

The simplest closed first-order formulas are identities,
which define varieties. Most of the results in~\cite{Ano02}
concern the problem of deciding whether a given black-box
$\Sigma$-algebra belongs to a fixed finitely based variety
of $\Sigma$-algebras. Although in~\cite{Ano02} it is assumed
that each element of a black-box $\Sigma$-algebra is
represented by a unique bit string (and hence equality
testing can be performed without access to an oracle), many
of the results in that paper are valid without this
assumption. We state two of these results here. First, every
variety of nilpotent groups (resp., Lie rings or Lie
$R$-algebras) is deterministically G-decidable in the class
of all groups (resp., Lie rings or Lie $R$-algebras). Here,
$R$ is an arbitrary finitely generated commutative
associative ring with~$1$. This result follows
from~\cite[Corollary~5.1]{Ano02} because for classes of
groups, rings, and $R$-algebras, strong G-decidability in
the sense of~\cite{Ano02} implies deterministic
G-decidability in the above sense. Second, suppose $\cl
V_1,\dots,\cl V_k$ are varieties of groups such that for all
$i\in\{1,\dots,k\}$, we have $\cl V_i\subseteq\cl A\cl
N_{c_i}\cap\cl N_{d_i}\cl A$, where $c_i$ and $d_i$ are some
nonnegative integers. Then $\cl V_1\dots\cl V_k$ is
G-decidable in the class of all groups
(see~\cite[Proposition~5.2]{Ano02}).

\subsection{Our contributions}
\label{ss:ourcontrib}

Suppose $\Omega$ is a fixed finite set of finitary operation
symbols. An $\Omega$-expanded group is a group, written
additively but not necessarily abelian, with an
$\Omega$-algebra structure. This group (without the
$\Omega$-algebra structure) is called the additive group of
the $\Omega$-expanded group. We assume that the symbols for
addition ($+$), additive inverse ($-$), and zero ($0$) do
not belong to~$\Omega$. See also
Definition~\ref{d:Omegaexpgranditsaddgr}. Furthermore, an
$\Omega$-expanded group $H$ is said to be distributive if
all the fundamental operations of $H$ associated with
nonnullary operation symbols in $\Omega$ are distributive
over addition. See also Definition~\ref{d:disOmegaexpgr}.

Let $c$ be a constant greater than~$1$. We construct a
probabilistic polynomial-time black-box $\Omega$-expanded
group algorithm $B$ such that the following holds. Suppose
$\bb H=(H,\rho)$ is a black-box distributive
$\Omega$-expanded group, $n$ is the encoding length of $\bb
H$, and $s\in H^*$ is a generating system of~$H$. Then $B$
on input $(1^n,s)$ outputs a generating system of the
additive group of $H$ with probability at least $1-n/c^n$.
See Theorem~\ref{t:compgensysofaddgr} for a precise
statement; the algorithm $B$ is described just before this
theorem. This result, referred to as \emph{Result~1},
generalizes the result of Arvind, Das, and Mukhopadhyay
stated in Subsection~\ref{ss:prevwork}. The algorithm $B$ is
based on the same ideas as the closure algorithm described
just before Lemma~2.3.8 in~\cite{Ser03}. We use random
subsums (see Subsection~\ref{ss:groups}) to construct
generating systems with $O(n)$ elements for certain
subgroups of the additive group of~$H$. Note that the
algorithm of Arvind, Das, and Mukhopadhyay mentioned in
Subsection~\ref{ss:prevwork} is based on different ideas.

Result~1 has at least two interesting applications. First,
we obtain a probabilistic polynomial-time black-box
$\Omega$-expanded group algorithm $C$ such that the
following holds. Let $\bb H=(H,\rho)$, $n$, and $s$ be as in
Result~1 and let $t\in H^*$. Then $C$ on input $(1^n,s,t)$
outputs with probability at least $1-2n/c^n$ a generating
system of the additive group of the ideal in $H$ generated
by~$t$. See Theorem~\ref{t:compgensysofaddgrofid} for a
precise statement; the algorithm $C$ is described in the
proof of this theorem. This result generalizes the result of
Babai et al.\ on the existence of a probabilistic
polynomial-time black-box group algorithm for
problem~\ref{i:findgensysfornormcl} (see
Subsection~\ref{ss:prevwork}). Using the algorithm of Babai
et al.\ for problem~\ref{i:rednumbofelem} (see
Subsection~\ref{ss:prevwork}), we can reduce the number of
elements in a generating system for the additive group of
the above-mentioned ideal to $O(n)$ in polynomial time. In
fact, the algorithm $C$ on the above input just runs the
algorithm $B$ twice with different parameters.

Note that the algorithm $B$ (resp., $C$) always outputs a
system of elements of $H$ (resp., the ideal in $H$ generated
by $t$) provided that the input is correct. Therefore, for
any polynomially bounded function $\eta\colon\N\to\R_+$, the
success probability of these algorithms can be amplified to
at least $1-2^{-\eta(n)}$. As above, this can be done in a
standard way by running the algorithm multiple times and
concatenating the output systems.

Second, Result~1 together with the result of~\cite{Ano02}
stated in Lemma~\ref{l:decgensysofaddgr} yields the
following. Suppose $\cl V$ is a finitely based variety of
distributive $\Omega$-expanded groups such that the additive
group of any $\Omega$-expanded group in $\cl V$ is
nilpotent. Then $\cl V$ is G-decidable in the class of all
distributive $\Omega$-expanded groups. See
Theorem~\ref{t:decgensys} for a precise statement. In our
opinion, this result is the main motivation for proving
Result~1.

In particular, all our results stated in this subsection
hold when the distributive $\Omega$-expanded groups are
ordinary groups, rings, $R$-modules, or $R$-algebras, where
$R$ is a fixed finitely generated commutative associative
ring with~$1$. Rings and $R$-algebras may be here with or
without $1$, where $1$ is considered as a nullary
fundamental operation. See also Remark~\ref{r:ringmodalg}.

\subsection{Open questions and organization of the paper}

Denote by $\cl G$ the class of all groups. Here are two open
questions for future research:
\begin{itemize}
\item Does there exist a finitely based variety of groups
that is not G-decidable in~$\cl G$? Note that finitely based
varieties of semigroups that are not G-decidable in the
class of all semigroups do exist
(see~\cite[Section~7]{Ano02}). Also, there exists a
quasivariety that is defined within $\cl N_2$ by a single
quasi-identity and is not G-decidable in $\cl N_2$
(see~\cite[Subsection~1.3]{Ano02}).

\item Let $d$ be an integer such that $d\ge4$. Is the
Burnside variety $\cl B_d$ G-decidable in~$\cl G$? For any
$k\in\N\setminus\{0\}$, the \emph{Burnside variety} $\cl
B_k$ is defined within $\cl G$ by the identity $\Forall
x(x^k=1)$ (see~\cite[14.63]{Neu67}). Note that $\cl B_1$
consists only of trivial groups, $\cl B_2\subseteq\cl A$,
and $\cl B_3\subseteq\cl N_3$. In fact, the inclusion $\cl
B_3\subseteq\cl N_3$ was proved by Levi and van~der~Waerden
in~\cite{LvdW33}; see also~\cite[12.3.5 and item~(iv)
of~12.3.6]{Rob96}. Hence, by the first of two results stated
in the last paragraph of Subsection~\ref{ss:prevwork}, we
see that $\cl B_1$, $\cl B_2$, and $\cl B_3$ are
deterministically G-decidable in~$\cl G$.
\end{itemize}

The rest of this paper is organized as follows.
Section~\ref{s:prelim} contains notation, definitions, and
general results used in the paper. In
Section~\ref{s:mainres}, we prove our results stated in
Subsection~\ref{ss:ourcontrib}.

\section{Preliminaries}
\label{s:prelim}

In this paper, $\N$ denotes the set of all nonnegative
integers. Let $Y$ be a set and let $n\in\N$. We denote by
$Y^n$ the set of all (ordered) $n$-tuples of elements
from~$Y$. Furthermore, we put $Y^*=\bigcup_{i=0}^\infty
Y^i$. Note that $\emptyset^*$ consists only of the empty
tuple. We consider elements of $\B^*$ as bit strings. The
unary representation of $n$, i.e., the string of $n$ ones,
is denoted by~$1^n$. Note that $e$ always denotes the base
of natural logarithms. If $\rv z$ is a random variable
taking values in a finite or countably infinite set $Z$,
then we denote by $\supp\rv z$ the support of $\rv z$, i.e.,
the set $\{z\in Z\st\Prob{\rv z=z}\ne0\}$.

\subsection{Universal algebra}

For a detailed introduction to universal algebra, the reader
is referred to standard books, \eg,~\cite{Cohn81, Wech92}.

Throughout the paper, $\Sigma$ denotes a set of finitary
operation symbols. Each $\sigma\in\Sigma$ is assigned a
nonnegative integer called the \emph{arity} of $\sigma$ and
denoted by~$\ar\sigma$. A \emph{$\Sigma$-algebra} is a set
$H$ called the \emph{carrier} (or the \emph{underlying set})
together with a family $(\widehat\sigma\colon
H^{\ar\sigma}\to H\st\sigma\in\Sigma)$ of operations on $H$
called the \emph{fundamental operations}. We usually denote
a $\Sigma$-algebra and its carrier by the same symbol when
there is no confusion.

Let $H$ be a $\Sigma$-algebra. Then the fundamental
operation of $H$ associated with a symbol $\sigma\in\Sigma$
will be denoted by $\sigma^H$ or simply by~$\sigma$.
Furthermore, suppose $\Psi$ is a subset of~$\Sigma$. Then
the $\Psi$-algebra obtained from $H$ by omitting the
fundamental operations associated with symbols in
$\Sigma\setminus\Psi$ is said to be the \emph{$\Psi$-reduct}
of $H$ (or the \emph{reduct} of $H$ to $\Psi$). We denote
the $\Psi$-reduct of $H$ by~$H|_\Psi$. Also, let $S$ be a
set or system of elements of~$H$. Then $\alg S_\Psi$ denotes
the subalgebra of $H|_\Psi$ generated by~$S$. For simplicity
of notation, if $\alg S_\Psi=H|_\Psi$ (in other words, $S$
generates $H|_\Psi$), then we will write $\alg S_\Psi=H$.

Denote by $\Sigma_0$ the set of all nullary operation
symbols in~$\Sigma$. Whenever $\sigma\in\Sigma_0$, it is
common to write $\sigma$ instead of~$\sigma()$. We note that
if $\Sigma_0=\emptyset$, then a $\Sigma$-algebra may be
empty, as in the books~\cite{Cohn81, Wech92}. This is not
standard in universal algebra, but is convenient, at least
in our opinion. However, for this paper it does not really
matter much whether we allow empty $\Sigma$-algebras. This
is because our main results are concerned with distributive
expanded groups and any expanded group is nonempty.

Suppose $X=\{x_1,x_2,\dotsc\}$ is a countably infinite set
of symbols called variables. We assume that $x_1,x_2,\dotsc$
are distinct and are not in~$\Sigma$. Denote by
$\Tm_\infty(\Sigma)$ the set of all $\Sigma$-terms over~$X$.
This set is defined as the smallest set $T$ such that
$X\cup\Sigma_0\subseteq T$ and if
$\sigma\in\Sigma\setminus\Sigma_0$ and
$v_1,\dots,v_{\ar\sigma}\in T$, then the formal expression
$\sigma(v_1,\dots,v_{\ar\sigma})$ is in~$T$. Of course,
$\Tm_\infty(\Sigma)$ is a $\Sigma$-algebra under the natural
operations.

An \emph{identity} (or a \emph{law}) over $\Sigma$ is a
closed first-order formula of the form
$\Forall{x_1,\dots,x_m}(v=w)$, where $m\in\N$ and
$v,w\in\alg{\{x_1,\dots,x_m\}}_\Sigma$. We usually omit the
phrase ``over~$\Sigma$.'' We will write identities simply as
$v=w$, where $v,w\in\Tm_\infty(\Sigma)$, assuming that all
variables are universally quantified. A class $\cl V$ of
$\Sigma$-algebras is said to be a \emph{variety} if it can
be defined (within the class of all $\Sigma$-algebras) by a
set of identities. By the Birkhoff variety theorem (see,
\eg,~\cite[Chapter~IV, Theorem~3.1]{Cohn81}
or~\cite[Subsection~3.2.3, Theorem~21]{Wech92}), a class of
$\Sigma$-algebras is a variety if and only if it is closed
under taking subalgebras, homomorphic images, and direct
products.

In this paper, $\cl V$ denotes a variety of
$\Sigma$-algebras. Suppose $\cl W$ is a variety of
$\Sigma$-algebras such that $\cl V\subseteq\cl W$. Then $\cl
V$ is called \emph{finitely based within $\cl W$} if $\cl V$
can be defined within $\cl W$ by a finite set of identities.
Also, $\cl V$ is said to be \emph{finitely based} if it is
finitely based within the variety of all $\Sigma$-algebras.
It is evident that if $\cl W$ is finitely based, then $\cl
V$ is finitely based if and only if $\cl V$ is finitely
based within~$\cl W$.

\subsection{Expanded groups}
\label{ss:expgroups}

Throughout the paper, $\Omega$ denotes a set of finitary
operation symbols. Similarly to $\Sigma$, the arity of any
$\omega\in\Omega$ is denoted by $\ar\omega$ and
$\Omega_0=\{\omega\in\Omega\st\ar\omega=0\}$.

\begin{defi}[$\Omega$-expanded group and its additive
group]
\label{d:Omegaexpgranditsaddgr}
Let $\Gamma$ be the set of operation symbols consisting of
the standard symbols for addition ($+$), additive inverse
($-$), and zero~($0$). (In particular, the arities of $+$,
$-$, and $0$ are $2$, $1$, and $0$, respectively.) We assume
that $\Gamma\cap\Omega=\emptyset$ and
$\Sigma=\Gamma\cup\Omega$. Then a $\Sigma$-algebra $H$ is
called an \emph{$\Omega$-expanded group} (or simply
\emph{expanded group}) if $\Add H$ is a group under $+$ with
inverse operation $-$ and identity element~$0$. This group
is said to be the \emph{additive group} of $H$; we emphasize
that it is not necessarily abelian.
\end{defi}

When it comes to $\Omega$-expanded groups, we assume that
$\Sigma=\Gamma\cup\Omega$, where $\Gamma$ is as in
Definition~\ref{d:Omegaexpgranditsaddgr} (in particular,
$\Gamma\cap\Omega=\emptyset$). Suppose $\Sigma$ is such a
set and $a$ and $b$ are elements of a $\Sigma$-algebra.
Then, of course, we write $a+b$ and $-a$ instead of $+(a,b)$
and $-(a)$, respectively. We also use $a-b$ as shorthand for
$a+(-b)$.

The notion of $\Omega$-expanded group generalizes the notion
of $\Omega$-group introduced by Higgins
in~\cite[Section~2]{Hig56}. Namely, an $\Omega$-expanded
group is an \emph{$\Omega$-group} if and only if
$\Omega_0=\emptyset$ and $\omega(0,\dots,0)=0$ for all
$\omega\in\Omega$. Note that $\Omega$-groups are also known
as groups with multiple operators (see,
\eg,~\cite[Chapter~II, Section~2]{Cohn81}). We emphasize
that rings with $1$ can be considered as $\Omega$-expanded
groups, but not as $\Omega$-groups if $1$ is required to be
a nullary fundamental operation. The same holds for
$R$-algebras with $1$, where $R$ is a commutative
associative ring with~$1$.

\begin{defi}[distributive $\Omega$-expanded group]
\label{d:disOmegaexpgr}
An $\Omega$-expanded group $H$ is called \emph{distributive}
if for any $\omega\in\Omega\setminus\Omega_0$, the
fundamental operation $\omega^H$ is distributive over $+^H$,
i.e.,
\begin{align*}
&\omega(h_1,\dots,h_{i-1},a+b,h_{i+1},\dots,h_{\ar\omega})
\\&\qquad
=\omega(h_1,\dots,h_{i-1},a,h_{i+1},\dots,h_{\ar\omega})
+\omega(h_1,\dots,h_{i-1},b,h_{i+1},\dots,h_{\ar\omega})
\end{align*}
for all $i\in\{1,\dots,\ar\omega\}$ and
$a,b,h_1,\dots,h_{i-1},h_{i+1},\dots,h_{\ar\omega}\in H$.
\end{defi}

It is easy to see that if $\Omega_0=\emptyset$, then
distributive $\Omega$-expanded groups are the same as
distributive $\Omega$-groups. Distributive $\Omega$-groups
were defined by Higgins in~\cite[Section~4]{Hig56}.

Clearly, the classes of all $\Omega$-expanded groups and of
all distributive $\Omega$-expanded groups are varieties of
$\Sigma$-algebras. The variety of all $\Omega$-expanded
groups is finitely based. If $\Omega\setminus\Omega_0$ is
finite, then the variety of all distributive
$\Omega$-expanded groups is also finitely based.

Of course, groups are exactly $\emptyset$-expanded groups,
which are always distributive. It is easy to see that the
class of all rings (resp., rings with $1$) is a finitely
based variety of distributive $\{\cdot\}$-expanded groups
(resp., distributive $\{\cdot,1\}$-expanded groups). The
symbols $\cdot$ and $1$ have their standard meanings here.

Let $H$ be an $\Omega$-expanded group. Then any subalgebra
of $H$ as a $\Sigma$-algebra is called an
\emph{$\Omega$-subgroup} of~$H$. This term is commonly used
for subalgebras of $\Omega$-groups; we extend it to
subalgebras of $\Omega$-expanded groups.

\begin{defi}[ideal of an $\Omega$-expanded group]
A set $A\subseteq H$ is said to be an \emph{ideal} of $H$ if
\begin{itemize}
\item $A$ is a normal subgroup of $\Add H$ and

\item for any $\omega\in\Omega\setminus\Omega_0$,
$i\in\{1,\dots,\ar\omega\}$, $h_1,\dots,h_{\ar\omega}\in H$,
and $a\in A$, we have
\[
\omega(h_1,\dots,h_{i-1},a+h_i,h_{i+1},\dots,h_{\ar\omega})
-\omega(h_1,\dots,h_{i-1},h_i,h_{i+1},\dots,h_{\ar\omega})
\in A.
\]
\end{itemize}
\end{defi}

The \emph{additive group} of an ideal in $H$ is this ideal
considered as a subgroup of~$\Add H$. For any ideal $A$ of
$H$, let $\theta_A$ be the equivalence relation on $H$ whose
equivalence classes are the cosets of $A$ in~$\Add H$. Then
it is easy to see that $A\mapsto\theta_A$ is a one-to-one
correspondence between the set of all ideals of $H$ and the
set of all congruences on~$H$.

\begin{rem}
\label{r:idealsofdisOmegaexpgr}
Assume that the $\Omega$-expanded group $H$ is distributive.
Suppose $A$ is a subset of~$H$. Then it is evident that $A$
is an ideal of $H$ if and only if $A$ is a normal subgroup
of $\Add H$ and for all $\omega\in\Omega\setminus\Omega_0$
and $i\in\{1,\dots,\ar\omega\}$, we have
$\omega(H,\dots,H,A,H,\dots,H)\subseteq A$, where $A$ occurs
in the $i$th position.
\end{rem}

See~\cite[Section~3, Theorem~4A, and its corollary]{Hig56}
for a definition and characterization of ideals of
$\Omega$-groups.

\begin{rem}
\label{r:omegagpSaromegasubseteq}
Assume that $H$ is distributive. Then for any set
$S\subseteq H$ and any $\omega\in\Omega$, we have
$\omega(\gp
S^{\ar\omega})\subseteq\gp{\omega(S^{\ar\omega})}$. Indeed,
if $\ar\omega=0$, then this is trivial. Otherwise this
inclusion follows from the distributivity of $\omega^H$
over~$+^H$.
\end{rem}

Define the function $\tau$ from the set of all subsets of
$H$ to itself by
\begin{equation}
\label{e:tau}
\tau(S)
=\GP{S\cup\bigcup_{\omega\in\Omega}\omega(S^{\ar\omega})}
\quad\text{for any set }S\subseteq H.
\end{equation}
The $k$-fold composition of $\tau$ with itself is denoted by
$\tau^k$ ($k\in\N$). In particular, $\tau^0$ is the identity
function on the set of all subsets of~$H$. It is evident
that $S\subseteq\tau(S)$ for all $S\subseteq H$. Also, if
$S\subseteq T\subseteq H$, then $\tau(S)\subseteq\tau(T)$.

\begin{lem}
\label{l:tauSeqtaugpS}
Assume that $H$ is distributive. Then $\tau(S)=\tau(\gp S)$
for every set $S\subseteq H$.
\end{lem}
\begin{proof}
Let $S$ be a subset of~$H$. Since $S\subseteq\gp S$, we have
$\tau(S)\subseteq\tau(\gp S)$. Furthermore, by
Remark~\ref{r:omegagpSaromegasubseteq}, $\omega(\gp
S^{\ar\omega})\subseteq\gp{\omega(S^{\ar\omega})}$ for all
$\omega\in\Omega$. Hence,
\[
\gp S\cup\bigcup_{\omega\in\Omega}\omega(\gp S^{\ar\omega})
\subseteq\gp S\cup\bigcup_{\omega\in\Omega}
\gp{\omega(S^{\ar\omega})}\subseteq\GP{S\cup
\bigcup_{\omega\in\Omega}\omega(S^{\ar\omega})}=\tau(S).
\]
This implies that
\[
\tau(\gp S)=\GP{\gp S\cup\bigcup_{\omega\in\Omega}\omega(\gp
S^{\ar\omega})}\subseteq\tau(S)
\]
because $\tau(S)$ is a subgroup of~$\Add H$. Thus,
$\tau(S)=\tau(\gp S)$.
\end{proof}

\begin{lem}
\label{l:uniontauiSeqalgSSigma}
For any set $S\subseteq H$, we have
$\bigcup_{i=0}^\infty\tau^i(S)=\alg S_\Sigma$.
\end{lem}
\begin{proof}
Suppose $S$ is a subset of $H$ and
$U=\bigcup_{i=0}^\infty\tau^i(S)$. It is easy to see that if
$T\subseteq\alg S_\Sigma$, then $\tau(T)\subseteq\alg
S_\Sigma$. This implies by induction on $i$ that
$\tau^i(S)\subseteq\alg S_\Sigma$ for all $i\in\N$.
Therefore, we have $U\subseteq\alg S_\Sigma$.

Let $\sigma\in\Sigma$ and $h_1,\dots,h_{\ar\sigma}\in U$.
Since $\tau^0(S)\subseteq\tau^1(S)\subseteq\dotsb$, we can
choose a positive integer $k$ such that
$h_1,\dots,h_{\ar\sigma}\in\tau^k(S)$. Then
$\sigma(h_1,\dots,h_{\ar\sigma})\in\tau^k(S)\subseteq U$ if
$\sigma\in\Gamma$ (because $\tau^k(S)$ is a subgroup of
$\Add H$) and
$\sigma(h_1,\dots,h_{\ar\sigma})\in\tau^{k+1}(S)\subseteq U$
if $\sigma\in\Omega$. This shows that $U$ is an
$\Omega$-subgroup of~$H$. Also, $S=\tau^0(S)\subseteq U$.
Hence we have $\alg S_\Sigma\subseteq U$. Thus, $U=\alg
S_\Sigma$, as required.
\end{proof}

\begin{rem}
\label{r:Rmodulesrepr}
Suppose $R$ is a commutative associative ring with $1$
generated by a system $s=(s_1,\dots,s_m)$ as a ring with
$1$, where $m\in\N$. In this remark, we assume that $\Omega$
consists of $m$ unary operation symbols
$\omega_1,\dots,\omega_m$. We show how to represent
$R$-modules as distributive $\Omega$-expanded groups. All
$R$-modules are assumed to be unital. Moreover, the class of
all $R$-modules represented in such a way is a finitely
based variety of distributive $\Omega$-expanded groups.
Similarly, the same can be done for $R$-algebras and
$R$-algebras with $1$, where $1$ is considered as a nullary
fundamental operation.

Let $B$ be an $R$-module. Consider $B$ as a distributive
$\Omega$-expanded group in which $\omega_i(b)=s_ib$ for all
$i\in\{1,\dots,m\}$ and $b\in B$. The additive group of $B$
as an $\Omega$-expanded group (i.e., $\Add B$) is the same
as that of $B$ as an $R$-module. It is evident that
$\omega_1^B,\dots,\omega_m^B$ are commuting endomorphisms
of~$\Add B$.

Suppose $y_1,\dots,y_m$ are independent commuting variables
and let $y=(y_1,\dots,y_m)$. As usual, denote by $\Z[y]$ the
polynomial ring in variables $y_1,\dots,y_m$ over~$\Z$. Put
$Q=\{p\in\Z[y]\st p(s)=0\}$. Then $Q$ is an ideal of $\Z[y]$
and $p+Q\mapsto p(s)$, where $p\in\Z[y]$, is an isomorphism
of $\Z[y]/Q$ onto~$R$. Of course,
$q(\omega_1^B,\dots,\omega_m^B)(b)=q(s)b=0$ for any $q\in Q$
and $b\in B$, where $\omega_1^B,\dots,\omega_m^B$ are
considered as commuting elements of the endomorphism ring
of~$\Add B$. The Hilbert basis theorem implies the existence
of a finite set $T\subseteq Q$ such that $Q=\sum_{t\in
T}\Z[y]t$. Choose such a set~$T$. Thus, we see that $B$
(considered as an $\Omega$-expanded group) satisfies the
following identities:
\begin{enumerate}[label=(\roman*)]
\item\label{i:abelgroup} $(x_1+x_2)+x_3=x_1+(x_2+x_3)$,
$x_1+x_2=x_2+x_1$, $x_1+0=x_1$, $x_1-x_1=0$;

\item $\omega_i(x_1+x_2)=\omega_i(x_1)+\omega_i(x_2)$ for
all $i\in\{1,\dots,m\}$;

\item $\omega_j(\omega_i(x_1))=\omega_i(\omega_j(x_1))$ for
all $i\in\{1,\dots,m\}$ and $j\in\{i+1,\dots,m\}$;

\item\label{i:tomegax1eq0}
$t(\omega_1,\dots,\omega_m)(x_1)=0$ for all $t\in T$.
\end{enumerate}
In identities~\ref{i:tomegax1eq0}, we use terms of the form
$p(\omega_1,\dots,\omega_m)(v)$, where $p\in\Z[y]$ and
$v\in\Tm_\infty(\Sigma)$. To avoid ambiguity, we formally
define these terms.

Put $I=\{(i_1,\dots,i_d)\st d,i_1,\dots,i_d\in\N\cs1\le
i_1\le\dots\le i_d\le m\}$. Suppose
$v\in\Tm_\infty(\Sigma)$. Then for every
$j=(i_1,\dots,i_d)\in I$, we define
$\omega_j(v)=\omega_{i_d}(\dots(\omega_{i_1}(v))\dots)$.
Furthermore, choose a total order (\eg, lexicographic)
on~$I$. For any nonempty finite multisubset
$J=\{j_1,\dots,j_k\}$ of $I$, where $j_1\le\dots\le j_k$,
let
\[
\omega_J(v)=(\dots(\omega_{j_1}(v)+\omega_{j_2}(v))+\dotsb)
+\omega_{j_k}(v).
\]
(Note that addition in $\Tm_\infty(\Sigma)$ is neither
commutative nor associative.)

Suppose $p\in\Z[y]$. Represent $p$ as
\[
\left(\sum_{(i_1,\dots,i_d)\in J_+}y_{i_1}\dots
y_{i_d}\right)-\left(\sum_{(i_1,\dots,i_d)\in
J_-}y_{i_1}\dots y_{i_d}\right),
\]
where $J_+$ and $J_-$ are disjoint finite multisubsets of
$I$; these multisubsets are uniquely determined by~$p$. Then
we define
\[
p(\omega_1,\dots,\omega_m)(v)=\begin{cases}
\omega_{J_+}(v)-\omega_{J_-}(v)&\text{if
}J_+\ne\emptyset\text{ and }J_-\ne\emptyset,\\
\omega_{J_+}(v)&\text{if }J_+\ne\emptyset\text{ and
}J_-=\emptyset,\\-\omega_{J_-}(v)&\text{if
}J_+=\emptyset\text{ and }J_-\ne\emptyset,\\
0&\text{if }J_+=J_-=\emptyset.
\end{cases}
\]

Conversely, let $H$ be a $\Sigma$-algebra satisfying
identities~\ref{i:abelgroup}--\ref{i:tomegax1eq0}, where
$\Sigma$ is as in Definition~\ref{d:Omegaexpgranditsaddgr}.
Then $H$ is a distributive $\Omega$-expanded group with
abelian additive group. Moreover, it is easy to see that $H$
is an $R$-module with additive group $\Add H$ and action
defined by $p(s)h=p(\omega_1^H,\dots,\omega_m^H)(h)$, where
$p\in\Z[y]$ and $h\in H$. This action is well defined
because $\omega_1^H,\dots,\omega_m^H$ are commuting elements
of the endomorphism ring of $\Add H$ and
$q(\omega_1^H,\dots,\omega_m^H)(h)=0$ for all $q\in Q$ and
$h\in H$. The latter follows from
identities~\ref{i:tomegax1eq0}.
\end{rem}

\subsection{Groups}
\label{ss:groups}

For a detailed introduction to group theory, we refer the
reader to standard textbooks, \eg,~\cite{Rob96, Rot95}.

Suppose $G$ is a group. A set of subgroups of the group $G$
is said to be a \emph{chain} if it is nonempty and totally
ordered by inclusion. If $C$ is a finite chain of subgroups
of $G$, then $\acl C-1$ is called the \emph{length} of~$C$.

\begin{rem}
\label{r:lenofchofsubgr}
It is well known that if $G$ is finite, then any chain of
subgroups of $G$ has length at most~$\log_2\acl G$. This
follows from the fact that if $K$ is a proper subgroup of a
finite group $L$, then $\acl L\ge2\acl K$.
\end{rem}

Let $g=(g_1,\dots,g_m)$, where $m\in\N$ and
$g_1,\dots,g_m\in G$. Then a \emph{random subproduct} of $g$
is defined as $g_1^{\re b_1}\dots g_m^{\re b_m}$, where $\re
b_1,\dots,\re b_m$ are chosen independently from the uniform
distribution on~$\B$. When $G$ is written additively, we use
the term ``random subsum'' instead. Namely, a \emph{random
subsum} of $g$ is $\re b_1g_1+\dots+\re b_mg_m$, where $\re
b_1,\dots,\re b_m$ are as above. Random subproducts were
introduced by Babai, Luks, and Seress in the preliminary
version of~\cite{BLS97} (see Subsection~6.2 in the
preliminary version and Definition~5.12 in the journal
version of that work). Some applications of random
subproducts can be found in~\cite[Section~2]{BCFLS95}. See
also~\cite[Section~2.3]{Ser03} and~\cite[Section~2]{CF93}.

The inequality in the next lemma is equivalent to
inequality~(2.5) in~\cite{Ser03} (see the proof of
Lemma~2.3.4 in that book). The original result is due to
Babai et al.~\cite[Theorem~2.5 and its proof]{BCFLS95}.

\begin{lem}
\label{l:Prrandsubprdonotgen}
Assume that the group $G$ is finite and that any chain of
subgroups of $G$ has length at most $l\in\N$. Suppose $s\in
G^*$ is a generating system of~$G$. Also, let
$k\in\N\setminus\{0,1,2\}$ and let $r_1,\dots,r_{kl}$ be
independent random subproducts of~$s$. Then
\[
\Prob{\{r_1,\dots,r_{kl}\}\text{ does not generate }G}\le
e^{-(1-2/k)^2kl/4}.
\]
\end{lem}

\subsection{Black-box universal algebra model}
\label{ss:bbunivalgmodel}

In this subsection, we assume that $\Sigma$ is an arbitrary
fixed finite set of finitary operation symbols. This allows
us to avoid representation issues.

\begin{defi}[black-box $\Sigma$-algebra]
Suppose $H$ is a finite $\Sigma$-algebra and $\rho$ is a
function from a subset of $\B^n$ for some $n\in\N$ onto~$H$.
Then the pair $(H,\rho)$ is called a \emph{black-box
$\Sigma$-algebra}.
\end{defi}

Let $\bb H=(H,\rho)$ be a black-box $\Sigma$-algebra. Any
algebraic statement concerning $\bb H$ should be understood
as concerning~$H$. For example, if $\cl C$ is an abstract
class of $\Sigma$-algebras, then $\bb H\in\cl C$ means that
$H\in\cl C$. Suppose $E$ is the domain of~$\rho$. If
$H\ne\emptyset$, then $E\ne\emptyset$ and there exists a
unique $n\in\N$ such that $E\subseteq\B^n$; this nonnegative
integer $n$ is said to be the \emph{encoding length} of~$\bb
H$. If $H=\emptyset$, then $E=\emptyset\subseteq\B^n$ for
all $n\in\N$. Therefore, in this case, any $n\in\N$ can be
taken as an encoding length of~$\bb H$. It is evident that
if $n$ is an encoding length of $\bb H$, then $\acl
H\le2^n$.

For every $h\in H$, we denote by $\rep h\rho$ an arbitrary
preimage of $h$ under~$\rho$. Any of these preimages will be
used as a representation of $h$ for computational purposes.
A similar notation was used by Boneh and Lipton
in~\cite{BL96}.

\begin{defi}[$\Sigma$-oracle]
\label{d:Sigmaoracle}
An oracle $O$ is called a \emph{$\Sigma$-oracle} for $\bb H$
if it works as follows:
\begin{enumerate}[label=(\roman*)]
\item\label{i:Sigmaoracle:equality} Given any query of the
form $(=,\rep{h_1}\rho,\rep{h_2}\rho)$, where $h_1,h_2\in
H$, the oracle $O$ decides whether $h_1=h_2$ (\eg, returns
$1$ if $h_1=h_2$ and $0$ otherwise).

\item\label{i:Sigmaoracle:operations} Given any query of the
form $(\sigma,\rep{h_1}\rho,\dots,\rep{h_{\ar\sigma}}\rho)$,
where $\sigma\in\Sigma$ and $h_1,\dots,h_{\ar\sigma}\in H$,
the oracle $O$ returns
$\rep{\sigma(h_1,\dots,h_{\ar\sigma})}\rho$.
\end{enumerate}
On other queries, the behavior of the oracle $O$ may be
arbitrary.
\end{defi}

To avoid confusion, when it comes to black-box
$\Sigma$-algebra model, we assume that the symbol $=$ is not
in~$\Sigma$.

\begin{defi}[black-box $\Sigma$-algebra algorithm]
A (possibly probabilistic) algorithm $A$ is said to be a
\emph{black-box $\Sigma$-algebra algorithm} if, when $A$
performs a computation in an arbitrary black-box
$\Sigma$-algebra,
\begin{itemize}
\item $A$ has access to a $\Sigma$-oracle for this black-box
$\Sigma$-algebra and

\item each query made by $A$ to this $\Sigma$-oracle has the
form specified in either item~\ref{i:Sigmaoracle:equality}
or item~\ref{i:Sigmaoracle:operations} of
Definition~\ref{d:Sigmaoracle} (of course, the item
specifying the form may depend on the query).
\end{itemize}
\end{defi}

When it comes to black-box $\Sigma$-algebra algorithms and
their computation in $\bb H$, we write $h$ instead of $\rep
h\rho$ for any $h\in H$. In particular, a tuple
$(h_1,\dots,h_m)$, where $m\in\N$ and $h_1,\dots,h_m\in H$,
and a set $S\subseteq H$ should be understood in this
context as $(\rep{h_1}\rho,\dots,\rep{h_m}\rho)$ and a set
$W\subseteq E$ such that $\rho(W)=S$, respectively. (Recall
that $E$ is the domain of~$\rho$.) This allows us to
simplify the notation. Of course, there exists a
deterministic polynomial-time black-box $\Sigma$-algebra
algorithm that, given a set $W\subseteq E$ and access to a
$\Sigma$-oracle for $\bb H$, computes a set $V\subseteq W$
such that $\rho(V)=\rho(W)$ and $\rho$ is one-to-one on~$V$.

Let $A$ be a probabilistic black-box $\Sigma$-algebra
algorithm. Consider a computation of $A$ in the black-box
$\Sigma$-algebra~$\bb H$. We will write $A^O$ to indicate
that $A$ uses a $\Sigma$-oracle $O$ for~$\bb H$. The
algorithm $A$ usually takes as input $(1^n,s)$, where $n$ is
an encoding length of $\bb H$ and $s\in H^*$ is a generating
system of~$H$.

We include $1^n$ in the input to the algorithm $A$ for the
following reason. For example, suppose $\Omega$ is a
nonempty finite set of nullary operation symbols. Consider
the case where $H$ is an arbitrary finite $\Omega$-expanded
group generated by~$\emptyset$. It is easy to see that $\Add
H$ is generated by $\{\omega^H\st\omega\in\Omega\}$. Hence
$\Add H$ can be isomorphic to any finite group generated by
a set of at most $\acl\Omega$ elements. Assume that for some
function $\beta\colon\N\to\R_+$, the running time of the
algorithm $A$ is bounded by $\beta(l)$, where $l$ is the
input length. (In particular, this holds if $A$ runs in
polynomial time.) Also, let the input to $A$ be the empty
tuple. Then $A$ should run in constant time. Therefore, the
power of the model seems to be very limited. For example,
suppose $O$ is a $\Sigma$-oracle for $\bb H$ and let
$\omega\in\Omega$. Then $A$ can query the oracle $O$ for
$\omega^H$, but cannot read the entire answer to this query
when $n$ is sufficiently large. This is because any
representation of $\omega^H$, as well as of every element in
$H$, has length~$n$. Including $1^n$ in the input to the
algorithm $A$ allows us to avoid this issue. We note that
$1^n$ is redundant in the input to $A$ in the following
cases:
\begin{itemize}
\item When the generating system $s$ is assumed to be
nonempty. We made this assumption in~\cite{Ano02}.

\item When $\Sigma$ contains no nullary operation symbols.
In this case, $\emptyset$ generates the empty
$\Sigma$-algebra.

\item When we are considering $\Omega$-groups in the sense
of Higgins~\cite[Section~2]{Hig56}, where $\Omega$ is a
finite set of nonnullary finitary operation symbols. (See
also Subsection~\ref{ss:expgroups} for a definition of an
$\Omega$-group.) In this case, $\emptyset$ generates a
trivial (i.e., one-element) $\Omega$-group. Note that this
case includes the cases of groups, rings, near-rings,
$R$-modules, and $R$-algebras, where $R$ is a finitely
generated commutative associative ring with~$1$.

\item When we are considering $\Omega$-expanded groups,
where $\Omega$ is a finite set of finitary operation
symbols, and assuming that $s$ is a generating system
of~$\Add H$. In this case, if $\Add H$ is generated by
$\emptyset$, then $H=\{0\}$.
\end{itemize}

\section{Main results}
\label{s:mainres}

In this section, we assume that $\Omega$ is a fixed finite
set of finitary operation symbols. Furthermore, suppose $c$
is a constant greater than~$1$. Choose an integer $k\ge3$
such that $c\le e^{(1-2/k)^2k/4}$. Such an integer exists
because $e^{(1-2/i)^2i/4}\to+\infty$ as $i\to+\infty$.

Let $B$ be a probabilistic polynomial-time black-box
$\Omega$-expanded group algorithm such that the following
holds. Suppose $\bb H=(H,\rho)$ is a black-box distributive
$\Omega$-expanded group, $n$ is the encoding length of $\bb
H$, $s\in H^*$ is a generating system of $H$, and $O$ is a
$\Sigma$-oracle for~$\bb H$. Then $B^O$ on input $(1^n,s)$
proceeds as follows:
\begin{enumerate}
\item Set $s_0=s$.

\item For each $i\in\{1,\dots,n\}$ (in ascending order), do
the following:
\begin{itemize}
\item Choose $kn$ independent random subsums of~$s_{i-1}$;
let $R_i$ be the set of all these subsums.

\item Compute the set
$S_i=R_i\cup\bigcup_{\omega\in\Omega}\omega(R_i^{\ar\omega})$.

\item Let $s_i$ be a sequence of all elements of $S_i$ taken
in an arbitrary order.
\end{itemize}

\item Return~$s_n$.
\end{enumerate}
When we need to specify $\Omega$, we write $B_\Omega$
instead of~$B$.

\begin{thm}
\label{t:compgensysofaddgr}
Suppose $B$ is the algorithm described just before this
theorem. Let $\bb H=(H,\rho)$, $n$, $s$, and $O$ be as in
the description of~$B$. Then $\supp B^O(1^n,s)\subseteq H^*$
and
\begin{equation}
\label{e:PrgpBO1nseqH}
\Prob{\gp{B^O(1^n,s)}=H}\ge1-\frac n{c^n}.
\end{equation}
\end{thm}
\begin{proof}
We will use the notation introduced in the description of
the algorithm~$B$. It is evident that $\supp
B^O(1^n,s)\subseteq H^*$. Since $\acl H\le2^n$,
Remark~\ref{r:lenofchofsubgr} implies that any chain of
subgroups of $\Add H$ has length at most~$n$.

\begin{clm}
\label{cl:Prforalli}
We have
\begin{equation}
\label{e:Prforalli}
\Prob{\Forall{i\in\{1,\dots,n\}}(\gp{R_i}=\gp{s_{i-1}})}
\ge1-\frac n{c^n}.
\end{equation}
\end{clm}
\begin{proof}
We have already seen that any chain of subgroups of $\Add H$
has length at most~$n$. Therefore, by
Lemma~\ref{l:Prrandsubprdonotgen}, we have
$\Prob{\gp{R_i}\ne\gp{s_{i-1}}}\le e^{-(1-2/k)^2kn/4}\le
c^{-n}$ for every $i\in\{1,\dots,n\}$. Hence,
\[
\Prob{\Exists{i\in\{1,\dots,n\}}(\gp{R_i}\ne\gp{s_{i-1}})}
\le\frac n{c^n},
\]
which is equivalent to inequality~\eqref{e:Prforalli}.
\end{proof}

In the rest of the proof, we use the function $\tau$ from
the set of all subsets of $H$ to itself defined
by~\eqref{e:tau}.

\begin{clm}
Assume that $\gp{R_i}=\gp{s_{i-1}}$ for all
$i\in\{1,\dots,n\}$. Then $\gp{s_j}=\tau^j(\gp s)$ for every
$j\in\{0,\dots,n\}$.
\end{clm}
\begin{proof}
We proceed by induction on~$j$. Since $s_0=s$, we have
$\gp{s_0}=\gp s=\tau^0(\gp s)$. If $j\in\{1,\dots,n\}$ and
$\gp{s_{j-1}}=\tau^{j-1}(\gp s)$, then
\[
\gp{s_j}=\tau(R_j)=\tau(\gp{R_j})=\tau(\gp{s_{j-1}})
=\tau(\tau^{j-1}(\gp s))=\tau^j(\gp s),
\]
where the second equality follows from
Lemma~\ref{l:tauSeqtaugpS}. Thus, the claim holds.
\end{proof}

\begin{clm}
\label{cl:taungpseqH}
We have $\tau^n(\gp s)=H$.
\end{clm}
\begin{proof}
It is evident that
\[
\gp s=\tau^0(\gp s)\subseteq\tau^1(\gp
s)\subseteq\dots\subseteq H.
\]
Let $l$ be the smallest nonnegative integer such that
$\tau^l(\gp s)=\tau^{l+1}(\gp s)$ (it exists because $H$ is
finite). Then by induction on $i$, $\tau^i(\gp s)=\tau^l(\gp
s)$ for all $i\ge l+1$. By
Lemma~\ref{l:uniontauiSeqalgSSigma}, $\tau^l(\gp
s)=\tau^{l+1}(\gp s)=\dots=\alg{\gp s}_\Sigma=H$. Since any
chain of subgroups of $\Add H$ has length at most $n$, we
have $l\le n$. Thus, $\tau^n(\gp s)=H$.
\end{proof}

Inequality~\eqref{e:PrgpBO1nseqH} follows immediately from
Claims~\ref{cl:Prforalli}--\ref{cl:taungpseqH}. Indeed, if
$\gp{R_i}=\gp{s_{i-1}}$ for all $i\in\{1,\dots,n\}$, which
holds with probability at least $1-n/c^n$, then
$\gp{B^O(1^n,s)}=\gp{s_n}=\tau^n(\gp s)=H$.
\end{proof}

For each $i\in\N$, $\omega\in\Omega\setminus\Omega_0$,
$j\in\{1,\dots,\ar\omega\}$, and
$d\in(\N\setminus\{0\})^{(\ar\omega)-1}$, let $\chi_i$ and
$\psi_{\omega,j,d}$ be unary operation symbols. We assume
that all these symbols are distinct and do not belong
to~$\Sigma$. For every $m\in\N$, put
\[
\Phi(m)=\{\chi_1,\dots,\chi_m\}\cup\{\psi_{\omega,j,d}\st
\omega\in\Omega\setminus\Omega_0\cs
j\in\{1,\dots,\ar\omega\}\cs
d\in\{1,\dots,m\}^{(\ar\omega)-1}\}.
\]
Note that $\acl{\Phi(m)}$ is bounded by a polynomial in~$m$.
Suppose $H$ is an $\Omega$-expanded group and
$g=(g_1,\dots,g_m)$, where $m\in\N$ and $g_1,\dots,g_m\in
H$. Let $H(g)$ be the $\Phi(m)$-expanded group such that
\begin{itemize}
\item $\Add{H(g)}=\Add H$ (in particular, $H(g)$ and $H$
have the same carrier),

\item $\chi_i^{H(g)}(h)=-g_i+h+g_i$ for all
$i\in\{1,\dots,m\}$ and $h\in H$, and

\item $\psi_{\omega,j,(d_1,\dots,d_{j-1},d_{j+1},\dots,
d_{\ar\omega})}^{H(g)}(h)=\omega(g_{d_1},\dots,g_{d_{j-1}},
h,g_{d_{j+1}},\dots,g_{d_{\ar\omega}})$ for all
$\omega\in\Omega\setminus\Omega_0$,
$j\in\{1,\dots,\ar\omega\}$,
$d_1,\dots,d_{j-1},d_{j+1},\dots,d_{\ar\omega}\in\{1,\dots,m\}$,
and $h\in H$.
\end{itemize}
It is evident that if $H$ is distributive, then $H(g)$ is
also distributive.

\begin{lem}
\label{l:idofHandPhimsubgrofHg}
Assume that $H$ is distributive. Then any ideal of $H$ is a
$\Phi(m)$-subgroup of $H(g)$. If $\alg g_{\{+,0\}}=H$, then,
conversely, every $\Phi(m)$-subgroup of $H(g)$ is an ideal
of~$H$.
\end{lem}
\begin{proof}
The former statement follows from
Remark~\ref{r:idealsofdisOmegaexpgr}. Now assume that $\alg
g_{\{+,0\}}=H$. Suppose $A$ is a $\Phi(m)$-subgroup
of~$H(g)$. Then $A$ is a subgroup of~$\Add H$. Moreover,
this subgroup is normal in $\Add H$ because if $f\in H$ and
$a\in A$, then $f=g_{i_1}+\dots+g_{i_k}$ for some
$i_1,\dots,i_k\in\{1,\dots,m\}$ ($k\in\N$) and hence
$-f+a+f=\chi_{i_k}(\dots(\chi_{i_1}(a))\dots)\in A$. Also,
since $H$ is distributive, we have
\[
\omega(h_1,\dots,h_{j-1},a,h_{j+1},\dots,h_{\ar\omega})
\in\alg{\{\psi_{\omega,j,d}(a)\st
d\in\{1,\dots,m\}^{(\ar\omega)-1}\}}_{\{+,0\}}\subseteq A
\]
for all $\omega\in\Omega\setminus\Omega_0$,
$j\in\{1,\dots,\ar\omega\}$,
$h_1,\dots,h_{j-1},h_{j+1},\dots,h_{\ar\omega}\in H$, and
$a\in A$. Therefore, by
Remark~\ref{r:idealsofdisOmegaexpgr}, $A$ is an ideal
of~$H$. Thus, the latter statement of the lemma holds.
\end{proof}

\begin{thm}
\label{t:compgensysofaddgrofid}
There exists a probabilistic polynomial-time black-box
$\Omega$-expanded group algorithm $C$ such that the
following holds. Suppose $\bb H=(H,\rho)$ is a black-box
distributive $\Omega$-expanded group, $n$ is the encoding
length of $\bb H$, $s\in H^*$ is a generating system of $H$,
and $O$ is a $\Sigma$-oracle for~$\bb H$. Also, let $t\in
H^*$ and let $Q$ be the ideal of $H$ generated by~$t$. Then
$\supp C^O(1^n,s,t)\subseteq Q^*$ and
\begin{equation}
\label{e:PrgpCO1nsteqQ}
\Prob{\gp{C^O(1^n,s,t)}=Q}\ge1-\frac{2n}{c^n}.
\end{equation}
\end{thm}
\begin{proof}
Suppose $C$ is a probabilistic polynomial-time black-box
$\Omega$-expanded group algorithm such that the following
holds. Let $\bb H=(H,\rho)$, $n$, $s$, $O$, $t$, and $Q$ be
as in the statement of the theorem. Then $C^O$ on input
$(1^n,s,t)$ proceeds as follows:
\begin{enumerate}
\item Run $B_\Omega^O$ on input $(1^n,s)$; let $g$ be the
output. (By Theorem~\ref{t:compgensysofaddgr} for
$B_\Omega$, we have $g\in H^m$ for some unique $m\in\N$.)

\item Run $B_{\Phi(m)}^{O_g}$ on input $(1^n,t)$ and return
the output, where $O_g$ is the $(\Gamma\cup\Phi(m))$-oracle
for $(H(g),\rho)$ naturally implemented by the oracle~$O$.
(It is easy to see that this can be performed in time
polynomial in $m$ and the length of $(1^n,t)$.)
\end{enumerate}

In the rest of the proof, we use the notation introduced in
the description of the algorithm~$C$. Suppose $T$ is the
$\Phi(m)$-subgroup of $H(g)$ generated by~$t$. Then by
Theorem~\ref{t:compgensysofaddgr} for $B_{\Phi(m)}$, we have
$\supp C^O(1^n,s,t)=\supp B_{\Phi(m)}^{O_g}(1^n,t)\subseteq
T^*$. However, Lemma~\ref{l:idofHandPhimsubgrofHg} implies
that $T\subseteq Q$. Hence $\supp C^O(1^n,s,t)\subseteq
Q^*$.

Assume that $\alg g_\Gamma=H$. Then $\alg g_{\{+,0\}}=H$
because $\Add H$ is a finite group. By
Lemma~\ref{l:idofHandPhimsubgrofHg}, ideals of $H$ are the
same as $\Phi(m)$-subgroups of~$H(g)$. Therefore, we have
$T=Q$. This shows that $\gp{B_{\Phi(m)}^{O_g}(1^n,t)}=T$ if
and only if $\gp{C^O(1^n,s,t)}=Q$. Thus, if
$\gp{C^O(1^n,s,t)}\ne Q$, then $\alg g_\Gamma\ne H$ or
$\gp{B_{\Phi(m)}^{O_g}(1^n,t)}\ne T$. Using this and
Theorem~\ref{t:compgensysofaddgr} for $B_\Omega$ and
$B_{\Phi(m)}$, we see that
\[
\Prob{\gp{C^O(1^n,s,t)}\ne
Q}\le\Prob{\gp{B_\Omega^O(1^n,s)}\ne
H}+\Prob{\gp{B_{\Phi(m)}^{O_g}(1^n,t)}\ne T}\le\frac
n{c^n}+\frac n{c^n}=\frac{2n}{c^n},
\]
which proves inequality~\eqref{e:PrgpCO1nsteqQ}.
\end{proof}

The next lemma follows from~\cite[Theorem~5.1]{Ano02}
because strong A-decidability implies deterministic
A-decidability. See~\cite[Section~3]{Ano02} for definitions
of these notions. We note that in~\cite{Ano02}, distributive
expanded groups are called multirings. The term
``multiring'' now usually means, loosely speaking, a ring
with multivalued addition (see,
\eg,~\cite[Definition~2.1]{Mar06}).

\begin{lem}
\label{l:decgensysofaddgr}
Let $\cl V$ be a finitely based variety of distributive
$\Omega$-expanded groups such that the additive group of any
$\Omega$-expanded group in $\cl V$ is nilpotent. Then there
exists a deterministic polynomial-time black-box
$\Omega$-expanded group algorithm $A$ such that the
following holds. Suppose $\bb H=(H,\rho)$ is a black-box
distributive $\Omega$-expanded group, $n$ is the encoding
length of $\bb H$, $g\in H^*$ is a generating system of
$\Add H$, and $O$ is a $\Sigma$-oracle for~$\bb H$. Then
$A^O$ on input $(1^n,g)$ decides whether $H\in\cl V$.
\end{lem}

\begin{thm}
\label{t:decgensys}
Let $\cl V$ be a finitely based variety of distributive
$\Omega$-expanded groups such that the additive group of any
$\Omega$-expanded group in $\cl V$ is nilpotent. Then there
exists a probabilistic polynomial-time black-box
$\Omega$-expanded group algorithm $D$ such that the
following holds. Suppose $\bb H=(H,\rho)$ is a black-box
distributive $\Omega$-expanded group, $n$ is the encoding
length of $\bb H$, $s\in H^*$ is a generating system of $H$,
and $O$ is a $\Sigma$-oracle for~$\bb H$. Then $D^O$ on
input $(1^n,s)$ correctly decides whether $H\in\cl V$ with
probability at least $1-n/c^n$.
\end{thm}
\begin{proof}
Choose a deterministic polynomial-time black-box
$\Omega$-expanded group algorithm $A$ according to
Lemma~\ref{l:decgensysofaddgr}. Let $D$ be a probabilistic
polynomial-time black-box $\Omega$-expanded group algorithm
such that if $\bb H=(H,\rho)$, $n$, $s$, and $O$ are as in
the statement of the theorem, then
$D^O(1^n,s)=A^O(1^n,B^O(1^n,s))$. By
Theorem~\ref{t:compgensysofaddgr} and
Lemma~\ref{l:decgensysofaddgr}, the algorithm $D$ satisfies
the desired condition.
\end{proof}

\begin{rem}
\label{r:ringmodalg}
Let $R$ be a fixed finitely generated commutative
associative ring with~$1$. In particular, all the theorems
in this section hold when the distributive $\Omega$-expanded
groups are ordinary groups, rings, rings with $1$,
$R$-modules, $R$-algebras, or $R$-algebras with~$1$. Here,
$1$ is considered as a nullary fundamental operation.
Furthermore, $R$-modules, $R$-algebras, and $R$-algebras
with $1$ are assumed to be represented as distributive
$\Omega$-expanded groups, where $\Omega$ is a suitable
finite set of finitary operation symbols. For $R$-modules,
this representation is described in
Remark~\ref{r:Rmodulesrepr}; for $R$-algebras and
$R$-algebras with $1$, the representation can be constructed
similarly. Note that $R$-modules are assumed to be unital,
as in Remark~\ref{r:Rmodulesrepr}.
\end{rem}

\bibliographystyle{plain}
\bibliography{ptime_algorithms_for_bb_distr_exp_groups}
\end{document}